# Idealized Social Dynamics

# In Bayesian Space of Assessments

Michael Zhdanov

Gainesville, Florida, U.S.A.

The purpose of this paper is to present an idealized hypotheses-drawn model describing dynamics of variability in business and other social areas. A new construct is introduced, the Bayesian space of assessments, to consider changes in positions of individuals and intrinsically interrelated entities called social bodies. The bodies' spatial coordinates are their assessments in any required for a specific purpose professional and/or ethical dimensions. A concept of market power, introduced originally in entrepreneurship and technology commercialization (Danov, Smith, and Mitchell, 2003), is extended to become a "driving force" applicable to other social processes. The model describes interrelations between driving forces acting on social bodies and changes in positions of these bodies in the Bayesian space of assessments. Implications are discussed for some business and political oscillations.

Matters of social, business, or group dynamics have been very popular for quite a long time in business and other social research and educational areas (e.g., Barile, Pellicano, and Polese, 2018; Forsyth, 2017). Tsoukas (2019) even considers "that everything is dynamic" and that "the dynamic nature is what is worth studying" (Vanderkerckhove, 2020). Leo Tolstoy (2010, originally 1867) "was on to behavioral ethics before there was such a thing as behavioral ethics" (Michaelson, 2022). In the same novel "War and Peace", Tolstoy introduced a term "people motions" meaning actions that may or may not be associated with their real motions in ordinary physical space and stated that the nature of people motions had yet been unknown. This is also seen from a comprehensive sociological review by Sorokin (1928), grasping, at least, a 60-year time span.

Understanding the nature of people motions refers now to the area of quantitative statistical or, in essence, high quality probabilistic empirical research and its stochastic analysis (e.g., Tuma and Hannan, 1984; Hegselmann and Krause, 2002; Toscani, 2006; Imai, 2017; Bellomo, Ha, and Outada, 2020; Zajdela et. al., 2022). Standing somewhat aside from the above is an important contribution to the discussion, the so-called "constructal theory" (e.g., Bejan and Merks, 2007). Like some of the above-mentioned theories, it is drawn on optimization principles and is considered common to both social and physical subjects.

In spite the above mentioned research directions are often called social dynamics, they rather represent kinematics, i.e., variability or change as opposed to statics, because they don't answer the foremost question of any dynamics of what sources of motions are, i.e., their driving forces, and how they affect changes in states of people affairs in eyes of other people, i.e., their motions in a space of assessments. The situation is, in a sense, somewhat like one was with physical sciences before Isaak Newton discovered interrelations between forces acting on physical bodies and their "states of affairs", spatial coordinates in the regular three-dimensional Euclidean space.

Bridging this gap or discovering functional interrelations between forces acting on social bodies and their states in the introduced here Bayesian space of assessments is the major purpose of the presented paper because that would have profound effects on further developments in social sciences, including business ones. This task seems to be feasible owing to just one frequently encountered phenomenon – social oscillations. It is performed through extending entrepreneurship and marketing concepts of transaction triangle and market power into wider business and other social areas (see also Danov, 2024).

By definition (Mitchell, 2001, 2002), transaction triangle concept elements "the creative entity", "others", and "the work" are intrinsically interrelated. A major force driving changes in their relationships is conceptualized as the market power (Danov, Smith, and Mitchell, 2003). These two concepts, transaction triangle and market power, are of the same importance for exchange processes in other areas of social life, including the political one, where a transaction triangle element "creative entity" is, for example, an election candidate or party and the "work" is their professional and/or ethical qualities, action programs, intentions, etc. Therefore, market power may work as a force, driving, in particular, historic changes or people motions.

Based on this idea, a deterministic analytical model is introduced to describe effects of market power and other driving forces on states of social bodies in the Bayesian space of assessments to understand frequently encountered oscillatory social phenomena. Implications are discussed, including trends in fashion industry, fluctuations in popularity between competing political parties in democratic societies, social sustainability, and interdependence between a society's political and economic structures and its overall performance.

ASSUMPTIONS AND CONSTRUCTS

*Subjects*

Subjects of the study, interacting elements of a transaction triangle, can represent individuals, groups, organizations (e.g., Forsyth, 2017; Robbins and Judge, 2019), network or market webs (Tjosvold and Weicker, 1993; Ryans, More, Barclay, and Deutcher, 2000), relationship webs (Danov, Smith, and Mitchell, 2003), or "a system", a generalization of the network theory (e.g., Golinelli, 2010; Barile, Pellicano, and Polese, 2018).

Instead of the mentioned above constructs of an individual, group, organization, or network, concepts of a social body and social group are applied, where the social body is a single social individual or a series, including interrelated individuals and other, smaller, interrelated social bodies, possessing or managing together common assets and/or sharing common values or goals, while social group is a series of unrelated social individuals or social bodies possessing just some common features and/or properties. Therefore, the group is a set of independent social elements, which motion could be better described with a stochastic approach, while a deterministic one may well be a better fit for social bodies representing a set of interrelated elements. Social bodies are the only subject of the presented paper.

This concept of group differs from its commonly used meaning of a series of interrelated individuals. Moreover, individuals and organizations are united here to possess common features similar to some other studies (Sepinwall, 2019; Martin, 2022). Social individuals are ones who tend to communicate with each other either directly or indirectly; live or work together as related elements; and share some material and/or spiritual values, goals, etc. Examples are humans, wolves, dolphins, swans, ants, bees etc. The only research subject here is humans, while some relevant studies (e.g., Szuba, 2001) consider also other social animals, insects, etc.

Examples of human social bodies are the following: individual; family; manufacturing, service, government, or non-government organization; county; state; industry; country; and the whole of humanity. While industry is an example of a social body, market is an example of a social group including unrelated parts. Other examples of human social groups are children, teenagers, seniors; laborers, farmers, teachers, engineers, scientists; or just a set of people in all countries aged between, say, 20 and 30.

Most large social bodies are built of smaller ones. For example, families are composed of individuals; local communities, of families or their parts; counties, of local communities; states or provinces, of counties; countries, of states or provinces; and humanity, of countries. The latter ones include industries consisting of companies. They also include numerous religious, sports, culture associations and organizations. Humanity, in addition to countries, includes international alliances and organizations. It has been mentioned already that any social body may be part of or may consist of many other social bodies. The ones including other social bodies, consisting of more than one individual, are called here complex social bodies.

*Space*

"People motions" can now be defined as social bodies' actions due to changes in their rational and emotional or spiritual states. These actions may or may not be associated with the bodies' physical motions in the regular visible three-dimensional space, but they do affect and are affected by their positions in an emotional/spiritual space of explicit and implicit assessments. A concept of moral or spiritual space was first constructed by the so-called "social physicists" in the seventeenth century (Sorokin, 1928).

A concept of space of probabilistic assessments (in a sense like the Bayesian interpretation of probability), introduced and employed here, is somewhat similar to the well-known concept of ratings widely used in analyzing public-opinion polls in a variety of business and other social areas, including the so-called "quantitative social science" (Imai, 2017). In spite that every individual assessment is subjective, a large number of these assessments, turned into a single point in the space of assessments, results in objective positioning of the evaluated subject.

This space of assessments may be one-dimensional with just one, e.g., professional or practical attribute being a function of time for a social body, material or immaterial product or service. It may also be a two-dimensional space where one professional/practical and one ethical/aesthetic independent "orthogonal coordinates" may be used to track the motions in time. It may be, at last, a multi-dimensional space covering several important independent professional and/or ethical constructs such as levels of skills, mutual trust, commitment, relational bonds, etc. (e.g., Anderson and Weitz, 1989; Crosby, Evans, and Cowles, 1990; Morgan and Hunt, 1994; Smith and Barclay, 1997).

In spite that discussions are still on about this topic (Harper, 2019), it is assumed that people motions are due to their decisions, explicit or implicit, made mostly intuitively of their emotions, imaginations, and moods rather than of their cognitions, thoughts, and rational analysis, while it is known that emotions themselves make our thinking more rational (e.g., Robbins and Judge, 2019), and, obviously, cognitions, thoughts, and rational analysis may generate imaginations, emotions, and moods.

This does not mean that rational or analytical is not important, but it means that sensual prevails rational. Indeed, one cannot construct a building, boat, or rocket based on just emotions and moods. An accurate mathematical analysis of the project, including its financial component, is a necessary requirement, but it is assumed here that a decision on taking the action is still drawn more on emotions and moods, while the conducted analysis affects them indirectly. Sometimes application of a force, like sudden appearance of an enemy, may cause quickly growing sense of fear, forcing the body to initiate an intuitive, unconscious action such as disorderly retreat, in spite that rational analysis could result in necessity of the opposite behavior.

If this is the case, and sensual does prevail rational, then subjective individual assessments obtained from representative polls and turned, eventually, into more objective ones, play the same role in the space of assessments as the role of spatial Euclidean coordinates in displaying motions of material bodies. The space of assessments may be one or two or multi-dimensional space including assessments of most important professional/practical and/or ethical/aesthetic dimensions for social/physical bodies respectively.

Therefore, results of applying driving forces in social processes, in general, and market power in commercial business processes, in particular, are positioned here in a space, derived from intuitive personal assessments, and social bodies' actions or inactions are both causes and effects of changes in their positions in the space of assessments. The nature of "people motions" is, in essence, reduced to the nature of forces driving emotions and moods of affected social bodies.

It is seen from the above that, unlike dynamics of material bodies, social dynamics cannot be studied in a regular coordinate space, and well-visible "people motions", such as legal and illegal changes of governments, insurrections, wars, revolutions, etc., may or may not be related to any significant spatial movements of people masses. Therefore, physical people displacements are not a subject of the current research, while motions of their spiritual states are. To study them, a new construct is introduced here. It is called the Bayesian space of probabilistic assessments, or just the Bayesian space of assessments, and it is similar, to some extent, to the mentioned earlier concept of moral or spiritual space introduced in the seventeenth century.

To create a coordinate system in a Bayesian space of assessments, the respective field research would have to be focused on more sophisticated polls than ones requiring just "approve-disapprove" answers, often expected in studies of ratings. The interviewed individuals would have to be asked about assessing a social/physical body's professional/practical and/or ethical/aesthetic or other attributes in a range between -100 and +100 with zero ratings assigned to unknown bodies only.

The range of [-100; +100] is chosen because an assessment of 1 out of 100 may also mean a 1% performance rating. Individual assessments in this interval (-100; +100) can be treated as probabilities (in the Bayesian sense) that "the said is true". E.g., an assessment of 90 of a Presidential Candidate's professional qualities means that the candidate would succeed professionally at this capacity with the likelihood of 90%. A respective negative assessment of -90 would mean that the candidate won't succeed with the same likelihood. In addition, the participants may feel comfortable because a similar range of (0; +100) is widely used in the North American education to evaluate performance, and many survey participants would likely be more accustomed to it than to other assessment spans.

Total assessments, however, cannot be expressed in percentage points because an assessment relative to a complex social body, an evaluator, is defined here as a sum of assessments by all social bodies composing the considered complex one. Therefore, a particular assessment, e.g., relative to a country, may well exceed 100; the space of assessments varies virtually from - ∞ to + ∞; and other units are required to measure the assessments. If an assessment is obtained in the above-described way, the respective unit of measurements could be called, for example, "leo" after Leo Tolstoy, who was a major contributor to inspiring this study, and the most appropriate name for the entire area of variability is the Bayesian space.

Starting from a certain level of social bodies, included into an appraiser complex social body, where the considered individual or another evaluated social body becomes unknown, both horizontally and vertically, further assessments of the social body will stay constant because assessments at further horizontal or higher vertical levels are zeros, while assessments of public figures or organizations known within more social bodies will continue to change. A body unknown at the humanity level will have the same assessment relative to the humanity as its assessment relatively to a lower vertical level of its country. If it is unknown at the country level, its assessment will stay the same as at the state level, etc.

All assessments at all levels are made in a way similar to the way of making ratings, i.e., through the so-called "representative samples". However, an assessment of an individual or another social body relative to an assessor social body is defined here as a sum of assessments at all lower levels compared to the assessor, which the considered individual or another social body may be known in directly or by reference. Therefore, total assessments of public social bodies may differ from assessments of ordinary ones by orders of magnitude, while present day regular assessments, such as ratings or IQs, cannot differ so significantly from each other.

*Forces*

A result of assessing a social body by another one is a dot in the above introduced Bayesian space of assessments relative to the assessor. If the first body, which is under the evaluation, applies a "driving force" towards the second one, the evaluator, this would affect the dot position in the space of assessments due to the same force of opposite direction, applied to the first body itself. Because this always happens in practice, it is assumed that the same force of the opposite direction is automatically applied to the body under the evaluation, and this opposite influence moves the dot in the space of assessments relatively to the chosen evaluator.

One of the most impressive examples is the "Bolshevik" (Red Army) victory in the Russian Civil War of 1918-1920. Their enemy the White Army was composed of very well armed professional warriors: well-trained soldiers and famous fighters Don Cossacks lead by brilliant Russian military officers. All of them were very motivated, patriotic, experienced people, but they lost to the Red Army, composed mostly of untrained peasants and industrial laborers. The major Red Army's power was in support from most ordinary civilians.

Bolsheviks won those "elections", i.e., most Russian people, owing largely to their notorious slogan, "Land to peasants, factories to laborers, peace to soldiers". It acted as a force, similar, in a sense, to the described below market power, affecting emotions and moods of the "voters", i.e., the Russian population, overwhelming majority of which were uneducated and even illiterate. As it happened, the slogan was fake. None of these social groups received anything, but it worked as a force driving the country down into an uncompetitive, utopian political system called socialism, while the opposite force, applied to the communistic party, pulled them up on top of the power over one of the greatest empires of all times.

A concept of market power (Danov, Smith, and Mitchell, 2003) was introduced as attractiveness of an entrepreneurial offer to potential stakeholders to optimize the process of technology commercialization. It is generalized here to become applicable to any potential buyer in virtually any transaction process, including in nonmonetary exchanges, and defined as a noncoercive Driving Force ($DF$) proportional to expected Perceived Benefits ($PB$) and inversely proportional to expected Perceived Costs ($PC$), i.e.,

$$DF \sim PB/PC. \tag{1}$$

In a particular case of commercial businesses, $PC$ turns into the product or service price to the buyer and the monetary benefit to the seller, while in other exchange processes real costs are often unclear. A transaction may happen if $PB/PC > 1$, and it is less likely to happen in a case of $PB/PC < 1$. This generalized concept of force turns into the concept of market power in a particular case of technology commercialization.

Even through this introductory definition, the concept of driving force seems to be useful in analyzing some social transactions, such as the reelection in the Republic of Belarus of then President Lukashenko in 2020. Both from numerous uprisings, visible unrest in the country following the elections and reported even by pro-Lukashenko Russian media and from private discussions with relevant to the state ordinary people, it is seen that the announced election result of 80% in favor of Mr. Lukashenko is somewhat difficult to believe in.

However, the real support for the former President of Belarus or, in essence, for the autocracy versus democracy, was likely much higher than Ukrainian support for the then President Yanukovych against the so-called "Maidan Revolution" in 2014. While the rough, "sensual", implicit evaluation of the perceived benefits $PB$ by both nations was about the same for virtually the same culture and people in very similar conditions in 2020 and 2014 respectively, some social groups in Belarus evaluated perceived costs *PC* in 2020 much higher and more realistically than many Ukrainians did back in 2014, because Belarusians could already see the Ukrainian Maidan and after Maidan events and experiences.

The supporting autocracy social groups in Belarus included most vulnerable people such as seniors, families with small children, who did not want any unrest in the country or any decline in their own living standard, low even without that, in exchange for any benefits, such as inviolability of property, freedom, democracy, self-government, and self-respect, etc. This may well be a real reason behind the obvious failure of all attempts to resist the seemingly fake election results. It means that relatively many Belarusians in 2020, compared to Ukrainians in 2014, voted for Mr. Lukashenko against any significant changes, and *DF* < 1 for some social groups in the triangle of opposition-voters-democracy, while *DF* > 1 for them in the triangle government-voters-autocracy.

There are many types of forces acting on social bodies, including ones unrelated to "internal", social sources and coming from material nature such as weather, wildfires, hurricanes, affecting spiritual state of social animals and, in particular, humans. Other types originate in alive, animal environments but are not social, because they are typical for non-social animals too, such as iguanas, alligators, tigers, bears, etc. Examples are sexual forces, attraction of parents to small cubs and cubs, to parents.

One of the most important types of social forces is the power or force of social attraction, making social animals social. Its existence is proven by the very existence of social animals. Everyone has had a chance to enjoy feelings of interest for new people or representatives of new social bodies. This interest may cover various areas including professional affiliation, the newcomer's judgements on your personal values or on traditional values of your related social bodies, etc.

It is assumed here that these very social forces of attraction are so-called "central" forces acting along the line connecting the social bodies in the Bayesian space of assessments. If forces of social attraction are indeed central, their magnitude should be inversely proportional to squared distances between two interacting social bodies in the space of assessments (e.g., Kittel, Knight, and Ruderman, 1973). It means that

$$F_a \sim 1/r^2 = 1/(\sum_1^n (x_{2i} - x_{1i})^2), \qquad (2)$$

where $F_a$ is a force of attraction between two social bodies, directed from one body to another; r is the distance between these two social bodies; $x_{2i}$ and $x_{1i}$ are $i$-th coordinates of the bodies 2 and 1 respectively in an $n$-dimensional space of assessments; and the sign $\sum_1^n$ is a sum from $i=1$ to $i=n$.

Similar to material bodies, tending to minimize their potential energy in the field of the Earth's gravity to get into a stable equilibrium, an alive body also tends to take a position in the space of assessments where its state is relatively sustainable. This idea is the major background for the mentioned earlier constructal theory, which is fruitful in explaining where people motions may go. However, it does not answer questions why they do go this way; what their driving forces are; and how they affect people.

A major of introduced here forces driving alive bodies, including social ones, is the power or force of doubts. While social bodies are tending to find a state of relative sustainability, "peace of mind", multiple factors affect them in the real world, acting often in opposite directions, and the force of doubts within social bodies works as a restoring force making their states sustainable. It is, in essence, a manifestation of an old philosophic principle known as the Hegel unity and struggle of opposites. Under the assumption that all forces are caused by emotions, examples of effects of the restoring forces of doubts are struggles between fear and fury, greed and generosity, love and hate.

The force of doubts is called onwards the social elasticity force because it is opposite to the direction of displacement from a sustainable state. In the linear approximation, it may be considered proportional to the magnitude of the displacement. This means that the restoring elasticity force is equal to

$$\boldsymbol{F}_e = - k_e \cdot \boldsymbol{r}, \tag{3}$$

where $\boldsymbol{F}_e$ is a vector of the social elasticity force, $\boldsymbol{r}$ is a vector of the displacement from a sustainable state, $k_e$ is yet an unknown coefficient of social elasticity in the space of assessments, and the module $\boldsymbol{r}$ is a displacement $r$ from the sustainable state:

$$|\boldsymbol{r}| = r = (\sum_1^n (x_{1i} - x_{0i})^2)^{1/2}, \tag{4}$$

where $x_{0i}$ are coordinates of the sustainable state in the $n$-dimensional Bayesian space of assessments, and $x_{1i}$ are coordinates of a displacement from the sustainable state.

Almost any other force driving a social body out of its sustainable state and called a force of change, can be represented in the first, linear approximation as

$$\boldsymbol{F}_c = k_c \cdot \boldsymbol{r}, \tag{5}$$

where $\boldsymbol{r}$ is a vector of displacement of the body from its sustainable state, $k_c$ is a coefficient of change, and $\boldsymbol{F}_c$ is a driving force of change. This is true in a small vicinity of any point because almost any natural function can be expanded into a Taylor series in its small vicinity, and the very first and most significant term of this expansion is linear in the first approximation.

In the already mentioned example of the Russian Civil War, the force of change was the market power or driving force due to attractiveness of the quoted Bolshevics' slogan and their respective program to most civilians, and it was implemented through the Red Army, while the social elasticity force was caused by the White Army, trying to return Russia to a sustainable state it was before. Both those forces acted together in opposite directions. Therefore,

$$\boldsymbol{F}_{total} = \boldsymbol{F}_e + \boldsymbol{F}_c = (k_c - k_e) \cdot \boldsymbol{r} . \tag{6}$$

A position of the respective point in the space of assessments is sustainable if $k_c < k_e$, and it is non-sustainable if $k_c > k_e$. If $k_e = k_c$, $\boldsymbol{F}_{total} = 0$. It is assumed here that if no forces are applied to social bodies, their points in the space of assessments are either still or moving along straight lines with permanent speeds, zero accelerations and decelerations.

*Inertia*

An assumption made at the very end of the last paragraph is a well-known property of material bodies called inertia. If no force is applied to a body, or $\boldsymbol{F} = 0$, the body is in a state of rest or moves by inertia unchanged. This is a property of resistance to change. It is assumed here that this property is an attribute of social bodies too. All of them possess inertia exactly like material ones do.

This assumption results from our everyday experience. Children seem to be much less inertial than adults. They are like a blank sheet of paper where adults can write virtually anything. They simply ape adults to learn quickly, and they are good at it. In particular, they can easily learn a second language to speak freely with no accent. Adults, on average, are more resistant to change, critical, even skeptical, and don't favor new ideas the way children do.

It is assumed that a social body is in a state of rest at a certain moment of time $t_0$, and no forces are applied to change this position in the space of assessments. I.e., the body itself is inactive until the moment $t_1$, as well as another social body, evaluating it, is inactive too with respect to the considered social body. Then nothing will happen to the body's ratings or assessments, and they will stay at rest in the time segment $[t_0, t_1]$.

If a body moves with a certain speed $v_0$ at the moment of time $t_0$, and no forces are applied to it until the moment $t_1$, it will go on moving with the same speed $v_0$ in the time segment $[t_0, t_1]$. This means that the body's assessments will be changing with the same speed in spite of no action going on. While it is not known yet if any results of any surveys could confirm this invariability of speed, it seems reasonable that assessments or ratings of a social body can be changed with no forces applied to it or to its evaluator. This assumption may be true in the only case if no friction, damping or other attenuation forces are applied to the body from the environment as well.

An example is the very first Russian President Mikhail Gorbachev, who teared down the Iron Curtain. After he resigned, retired, and did nothing in politics, gave almost no interviews to media, his world ratings continued to grow. Moreover, they went on growing even after his death. Unlike in material nature, where dead bodies turn to nothing, social bodies in spiritual world may continue affect others for quite a long time owing to their life-time achievements. A well-known saying is true that a personality dies only when everyone forgets about it. One more applicable saying, confirming that ratings may change with no forces applied, is that "the great can be seen only from a distance".

A measure of social inertia is called social mass here. The greater a social body's mass, the greater its social inertia or resistance to change. It is assumed that the social mass is also a measure of the property of attraction between social bodies. This property was discussed here earlier, in the subsection *Forces*.

It seems that the social mass is a more sophisticated and difficult to measure attribute than the well-known material mass. In particular, the above-mentioned example of differences between children and adults shows that the social mass should be a function of intellectual development, which measure is more difficult to determine than, for example, a regular measure of intelligence called IQ. Quantitative differences in this new measure between kids and adults, between average, ordinary people, and human giants, such as Shakespeare, Tolstoy, Tchaikovsky, Renoir, etc., should be much, much greater than entire existing differences in IQ measures now (e.g., Goldstein, Allen, and Deluca, 2019). While applying to ones who lived in the past, inertia is only a measure of attraction to them and not a measure of their own resistance to change.

As both a measure of attraction and resistance to change, the mass should also be dependent on social bodies' physical abilities. They are sharp teeth and strong muscles for wolves, muscles and guns for humans, armed forces for countries, etc. Not by chance, there is a saying that "mind is not necessary if you got power". Finally, inertia should be a function of social bodies' assets, their economic and/or financial abilities. These properties could be measured, for example, in the value of assets and earnings for individuals and commercial organizations, assets and GDPs for counties, states and countries, etc.

Although all these three variables (intellectual, physical, and economic) are interdependent, their independence of each other in a certain range of their variability may often be observed in real societies. While even a unit to measure the social mass has not been introduced yet, it is qualitatively clear that the greater any one of the above three variables is, the greater the total mass of the respective social body.

Therefore, mass *m* of a social body is a function of three variables: $m = m(x, y, z)$, where *x* is an intellectual variable, *y* is a physical or military variable, and *z* is an economic variable of a social body. If $x = y = z = 0$, then $m = 0$, and m is not equal to 0 if any of these variables is not equal to 0.

Introduced this way mass could be measured, for example, in monetary units because it depends on three variables, and each of them can be expressed in monetary units for human social bodies. However, measuring a social mass is a special problem more difficult than measuring a mass of a physical body, and it is not a subject of this paper. Despite that, even its mere introduction results in valuable conclusions consistent with some well-known observations.

## MODEL

*Background*

Driving forces emerging from interactions within a transaction triangle make participating social bodies affected by either each other or work they produce according to transaction cognition theory, or by both (Mitchell 2001, 2002).

The forces change judgements, i.e., positions of the affected bodies in a Bayesian space of assessments, causing relationship changes and, in essence, social changes or "people motions". In addition to these forces, resulting from the bodies' interactions, there can be ones caused by independent "third" parties, i.e., bodies standing out of the considered transaction triangle, as well as by surrounding environments, e.g., unusual weather effects such as wildfires, hurricanes, etc., affecting feelings, moods, assessments, and, therefore, states of respective social bodies. This additional very close, but a somewhat wider ambit of social dynamics is not a subject of the current study.

The following three subsections, *Attraction*, *Dynamics*, and *Oscillations*, constitute the introduced here model of social dynamics. It is drawn on well-known variabilities of social bodies called oscillations.

Equations describing social oscillations are still unknown just because unlike with physical bodies, social ones move in a space of human assessments invisible to the naked eye. While results of these spiritual movements, such as new or broken business contracts, riots, revolutions, etc., are well seen, the motions themselves are yet invisible. They can be seen only in the already introduced here Bayesian space of assessments constructed through independent unanimous, confidential polls or through other yet less popular means such as direct measurements of brain responses, e.g., via polygraphs.

*Attraction*

Similar to gravitational forces in mechanics, there are, obviously, forces of attraction between all social bodies, including human ones. Otherwise, there would be no social animals themselves. These forces, according to an assumption made in the subsection *Forces*, relate to the so-called "central" forces, which are inversely proportional to squared distances between affected bodies in physical space and to squared differences in assessments between social bodies, described by the relation (2), in the Bayesian space of assessments.

This means that bodies with close ratings attract stronger to each other. Celebrities or so-called "stars" are usually more interested in each other than they are in ordinary people with low or no ratings. While ordinary social bodies are often more interested in each other than in public ones too, an ordinary individual is normally more attracted to a celebrity than the celebrity to the ordinary person. This happens because the social masses of these two types of bodies are very different. The phenomenon becomes clearer if the law of social attraction is similar to the law of gravitational attraction:

$$F_a = \gamma \cdot m_c \cdot m / r^2, \tag{7}$$

where $F_a$ is a module of the force of attraction between two social bodies, $r^2$ is the squared distance between them in the Bayesian space of assessments, $m_c$ and $m$ are social masses of the two bodies, being, in this case, masses of a larger body and much smaller one, a celebrity and an ordinary one, and $\gamma$ is a coefficient of attraction, independent on a particular body.

*Dynamics*

As it is mentioned in the introduction, the term "dynamics" is often misused in present-day social sciences because it stands for the "force" or "power" due to its meaning in the Greek language. However, it is often used to show just variability or change as opposed to statics, while in theoretical mechanics it does mean interdependence between forces applied to physical bodies and their motions in the regular three-dimensional space. In this, latter sense, a science of social dynamics does not yet exist as interdependence between driving forces applied to social bodies and their motions.

Following the introduction of the Bayesian space of assessments, where spiritual motions become visible, and discovery of driving forces capable of acting on social bodies, such as the force of social elasticity (3) and the force of social attraction (7), it becomes feasible to state major laws of social dynamics. These laws may well be applicable to interactions amongst transaction triangles, which, themselves, may be considered as complex social bodies, and to non-social, e.g., natural or environmental forces, acting on social bodies. However, the scope of the presented work does not cover the latter forces, and, for now, three major laws describing interactions within a transaction triangle are the following.

Firstly, it has already been assumed in the subsection *Inertia* that a social body would either stay at the same point or move with the same permanent speed and direction in the Bayesian space of assessments if no forces are applied. Secondly, if a driving force is applied to one of two interacting social bodies because of an ongoing relationship between them, then the same driving force of the opposite direction is applied to another party of the relationship. This is true, in particular, for the force of attraction (7).

And, at last, thirdly, any force applied to a social body causes an acceleration of this body in the Bayesian space of assessments equal to the magnitude of the applied force divided by the value of mass of this social body:

$$\boldsymbol{a} = \boldsymbol{F} / m \ . \tag{8}$$

Here $\boldsymbol{F}$ is a vector of force, applied to a social body, $\boldsymbol{a} = d\boldsymbol{v}/dt$ is a vector of its acceleration, $\boldsymbol{v} = d\boldsymbol{r}/dt$ is a vector of its speed, $\boldsymbol{r}$ is a vector of the social body's position in the Bayesian space of assessments, and $t$ is time.

It will be seen in the next subsection that in its so-called linear approximation, (8) turns into second order linear differential equations, which solutions describe small-amplitude oscillations and waves.

It has been well known for at least a hundred years that oscillations of social bodies do exist (e.g., Sorokin, 1928). Numerous modern observations also report their presence in various social areas related to stocks, markets, labor, fashion, politics, and history. Therefore, the law (8) must work in social dynamics just because all small-amplitude oscillations are mathematically described by second order linear differential equations.

Meanwhile, it will be seen here and now how well the law (8) explains a social phenomenon, mentioned earlier in the subsection *Attraction*. Combining (7) and (8), it is easy to obtain that

$$\gamma \cdot m_c \cdot m / r^2 = m_c \cdot a_c \qquad (9)$$

and

$$\gamma \cdot m_c \cdot m / r^2 = m \cdot a . \qquad (10)$$

Through dividing (9) by (10) obtain that

$$a_c / a = m / m_c . \qquad (11)$$

It is seen from the (11) that $a_c \ll a$ because $m_c \gg m$. If a larger body of mass $m_c$ is a celebrity, and a smaller one of mass $m$ is an ordinary individual, then, despite that attraction forces acting on both the celebrity and the ordinary body are equal in magnitude, the celebrity is much less attracted to the ordinary body because its social mass $m_c$ is much greater than the mass $m$ of the ordinary body.

*Oscillations*

Social oscillations caused by the struggle of opposites or social elasticity can be described in the first, linear approximation through combining the law of dynamics (8) with the expression (3) for the respective force of social elasticity to obtain a second order differential equation

$$m \cdot d^2 r / dt^2 + k_e \cdot r = 0, \qquad (12)$$

which is reduced to a scalar equation

$$d^2 x / dt^2 + (k_e/m) \cdot x = 0 \qquad (13)$$

for assessments depending on one-dimensional variables only. Here $x$ is one of these variables, e.g., a professional or ethical dimension of a "creative entity" or "others" or a practical or aesthetic dimension of "work".

The equation (13) is good to describe and interpret as oscillations customer "others" doubts in accepting an offer of an element "work" in a transaction triangle. It explains, for example, why fashions change because a fashion element is "work" by a fashion designer, which a social body "others" may or may not buy.

It is also helpful in understanding why there is a saying, "history repeats itself" with regards to small, "linear" historic oscillations such as legal fluctuations of power between two political parties in free, democratic societies. Employing full none-linearized equations (8) would enable one to model huge, memorable social oscillations as well as to forecast better upcoming large, non-linear social changes, including riots, revolutions, wars, virtually any large economic, political, or cultural events.

It is seen from the subsections *Subjects* and *Space* that large complex social bodies include smaller ones. Because forces of social attraction are central, the complex bodies are shaped as globes in a respective Bayesian space of assessments. For example, in a two-dimensional plane with axis $x$ being ethical and $y$, professional assessments, the globes turn into discs.

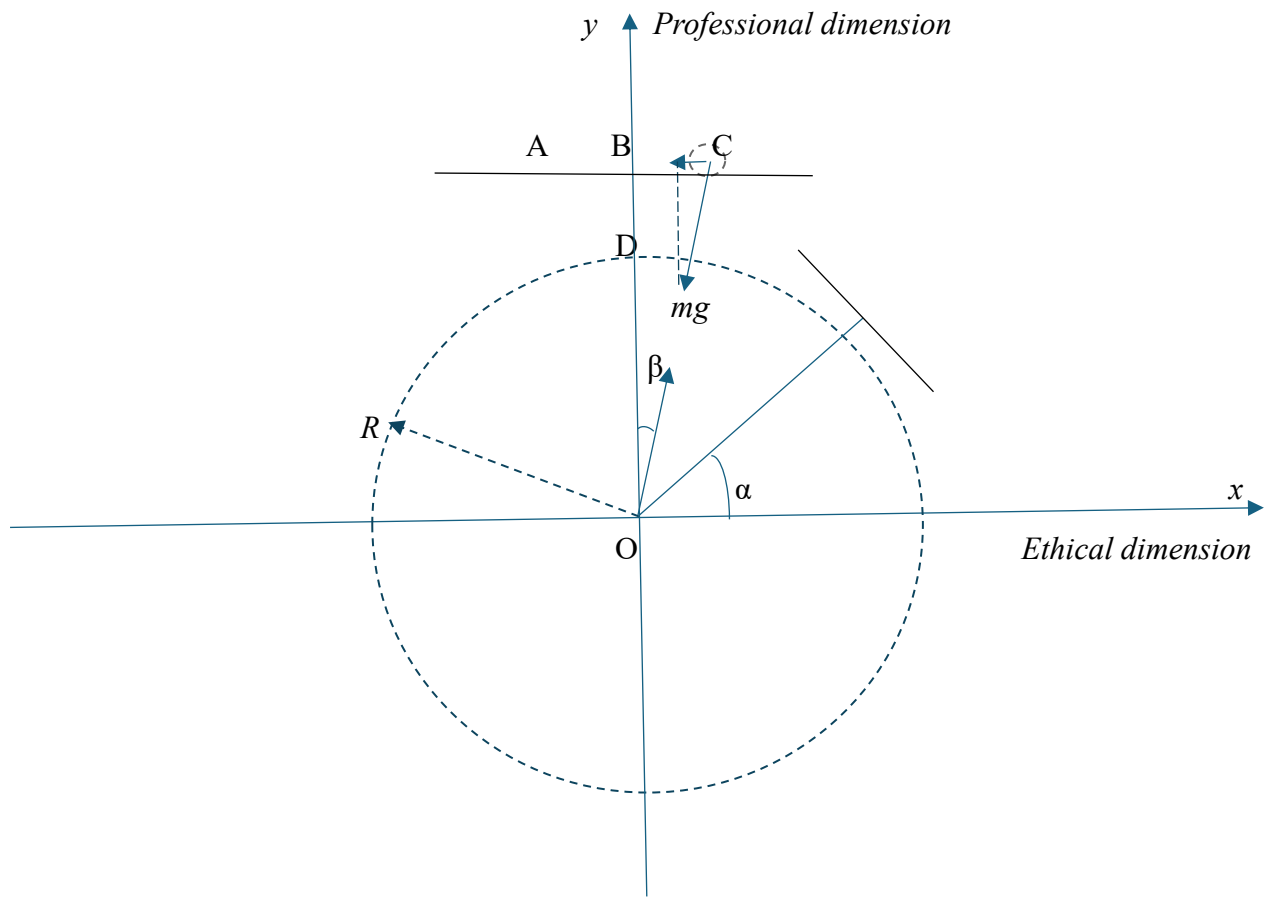

Figure 1. Social Attraction Oscillator

Fig. 1 describes the attraction between a social body C and a large, complex body O in a two-dimensional Bayesian space of assessments with axis $x$ being ethical and $y$, professional dimensions. A disc centered at the point O with coordinates (0,0) is a large, complex social body, e.g., a large country, such as the U.S., consisting of a numerous number of smaller social bodies including the one centered at the point C.

The body C is located close to the radius $R$ of the body O, and the forces of attraction between them are determined from the expression (7) with $r = R$, which is much greater than a module of assessment of an average, ordinary publicly unknown individual or another small social body such as a small and unknown private company, NGO, etc. It means that body C is a public body relative to the body O, i.e., it is known to many or all social bodies inside of the body O.

In considering so-called free oscillations, a particular example of the body C is a widely known in the country public individual such as a former President of the U.S., etc. The term "former" is not accidental. Only inactive public social bodies are considered in this example of free oscillations to avoid motions affected by any forces other than the one of attraction. It means that all achievements and results are in the past for the considered bodies, e.g., the considered former president. The area above the line $y$=OB is a locus of sets of the body's assessments relative to a large complex body O during a span starting from the end of the presidency, when its assessment is at the point C, until a moment $t$, when it is still inactive.

It means that owing to the body's previous achievements, there is a straight line of a minimal professional assessment, which the inactive body can get, and the line is working similar, in a sense, to a rigid slippery plane laying on the Earth surface to limit physical motions of a material body. Despite that any social body within the body O can give any assessment, and some very limited number of these assessments may fall below the line, it is assumed that eventual assessment obtained through summarizing evaluations by all social bodies will never fall below this line. If this line is a curve, then resulting oscillations or other motions are dependent on its curvature.

To go on, a concept of "ordinary public bodies" is defined relatively to a large complex social body, such as a large country, as the one "laying on the surface" of the body. I.e., their modules of assessments are all within a certain layer BD, which thickness $\Delta R$ is much smaller than the radius $R$, i.e., $\Delta R \ll R$. Then ordinary public bodies are ones laying in the intervals of $(R, R+\Delta R)$, while "ordinary bodies" are all within $(0, R)$, and "outstanding public bodies" are in the interval $(R +\Delta R, \infty)$. Only small motions of the ordinary public body C are considered, for which $tg\ \beta = BC/BO \ll 1$, i.e., the maximum displacement BC is much smaller than the radius $R$.

If the radius $R$ of a complex body is defined, for example, as just the total number of social bodies withing the considered complex body, then an ordinary public body known to all these bodies will be assigned the total evaluation of $R$, laying right on the surface of the considered complex body, if all its assessments by every one of the compounding bodies are positive and minimal, i.e., equal to one. If the complex body is, for example, a residential community, the radius is just the number of the community residents, their families and associated workers.

The whole situation is somewhat similar to the one with material bodies on Earth located within the troposphere BD, which thickness is about 10 km, while the Earth's radius is over 6000 km, and $BD/R = 1/600 \ll 1$. Due to the above given definition of the ordinary public body, the values of assessments relative to the social body O are about the same for all ordinary, not outstanding, public bodies, and they are close to $R$ for the example in Fig.1. This means that the mentioned ordinary public bodies are all in the layer $BD \ll R$. Therefore, distances are about the same between centers of the "ordinary" public bodies located within BD and the center of a large social body they all interact with. Then it is easy to obtain from (7) that the force acting on an ordinary public body is equal to

$$P = g \cdot m,  \quad (14)$$

where *P* is a force of social attraction acting on the body and directed towards the center of a the disc-like large complex social body O, such as a big and well-known country; $g = \gamma \cdot m_c / R^2$, with γ being a constant coefficient in the expression (7), $m_c$ is a mass of the large social body O, *R* is its radius, and *g* is a "free-fall" acceleration, which is about the same for all ordinary public bodies, whose assessments are all within the layer BD in the Fig. 1. This is a well-known fact for material bodies on Earth, and it is also true for ordinary public bodies interacting with a large social body of radius *R*, as far as *R* = *const* with high accuracy due to the inequality BD << *R* .

It is seen from Fig. 1 that a projection of the force of attraction *mg* to the axis *x* is equal to -*m*·*g*·Sin *ß* = -*m*·*g*·(*x*/*R*). The force is negative because its direction is opposite to the positive direction of the axes *x*. Then, due to the equation (8), $F = ma = m \, d^2x / dt^2 = - mg(x/R)$ and

$$d^2x / dt^2 + (g/R) x = 0 . \qquad (15)$$

The equations (15) and (13) are the same, except for the coefficients before *x*. Both of them describe small oscillations through their general solution

$$x = A \, Sin \, \omega t + B \, Cos \, \omega t , \qquad (16)$$

where $\omega = \sqrt{g/R}$ and $\omega = \sqrt{k_e/m}$ for the equations (15) and (13) respectively. Here ω is the circular frequency, ω = 2π /*T*, with *T* being the period of the oscillations, the regular frequency ν = 1/*T,* and, respectively, ω = 2πν.

The coefficients *A* and *B* are determined from initial conditions at the moment of time *t* = 0. It is seen, in particular, from the (16) that

$$B = x(0), \qquad (17)$$

and

$$A = dx/dt(0):\omega = v(0)/\omega, \qquad (18)$$

where $x(0)$ is a displacement of the body from the equilibrium at the moment of time $t=0$, and $v(0)$ is a speed of the body at the initial moment of time in the space of assessments. If, for example, at the moment of time $t=0$, the body is at the point C and $v(0)=0$, then $A=0$ and $B=BC=max(x)$. If the line OB is located at an angle α to the $x$-axis, as depicted in Fig. 1, then both professional and ethical assessments oscillate. The only professional attribute oscillates in a particular case of $α = 0$.

The equations (13) and (15) describe, through their solutions (16), only free oscillations, when a social body is taken out of its state of rest or straight-line, uniform motion. These oscillations are idealistic, perpetual motions. In reality, physical body oscillations gradually disappear in the material world due to always present damping caused by the body's interactions with the surrounding substances such as air, water, an underlying surface, etc.

The aforementioned is true for social motions too, in spite that the nature of the attenuation is different. In real life, every social body is part of several, not just one, transaction triangles, and decay of oscillations in one of them may be due to distraction of the considered social body to other triangles, causing gradual oblivion of the initial event generated the oscillations. The attenuation can be easily introduced into the equations (13) and (15). Their solutions, modified relations (16), are also known. The decay is not taken into account intentionally in this idealized approach to simplify and, therefore, clarify the presented introduction into the nature of social oscillations.

IMPLICATIONS FOR BUSINESS AND POLITICS

     While the existence of just one phenomenon – social oscillations – was a major motivation to pursue the presented research, its implications may well go far beyond this one of the most impressive elements in the big picture of the world we, humans, live in. In fact, the number of this model potential applications is hard to overestimate.

     An example of one of them – greater attractiveness of public social bodies to each other than to ordinary ones – has been interpreted in the subsection *Dynamics*. To start, however, just two more relatively simple, qualitative examples of applications are discussed onwards in the following subsections called *Fashion Fluctuations* and *Political Oscillations*. Interdependence between economic and political society structures, as well as a society's stability, are considered in the very last subsection *Social Sustainability*.

*Fashion Fluctuations*

     There is a well-known saying, "fashions always come back". To analyze this phenomenon in the light of the presented model, it is required to identify a relevant transaction triangle. The "creative entity" here is a social body called fashion designer, who produces a limited number of its new items of closing, "work", for a limited number of first customer "others".

     The "others" were, for example, nobilities, members of a royal court in the medieval France. They are so-called celebrities and other ordinary public social bodies in modern fashion-creator countries such as the U.S., France, Italy, etc. To make it simpler, the "work" is specified here as an element or parameter of closing. It may be, for example, the length of women skirts or

the width of trousers bottom. Everyone knows that the length goes up and down and the width increases and decreases over time. Averaging them over a long-time span, say over about a hundred years or so, provides with a normal level of an equilibrium state as well as with displacements relative to this average level. Then some of the displacement parameters can be evaluated. Their examples are the oscillations period $T$, its amplitude $A$, speed of change, etc.

These fluctuations are due to changing judgements by customer "others" about a "work" of a "creative entity", its length of skirts or trousers width. I.e., they are due to the power of doubts or social elasticity force acting on a social body representing the fashion designer "creative entity" through evaluation of its "work" by customer "others". The customer assessments $x$ of the length or the width are most likely neutral, i.e., $x=0$, at the state of equilibrium, and they should be reaching their maximum $x=A$ or minimum $x=-A$ at a state where the length or width are the shortest or longest ones.

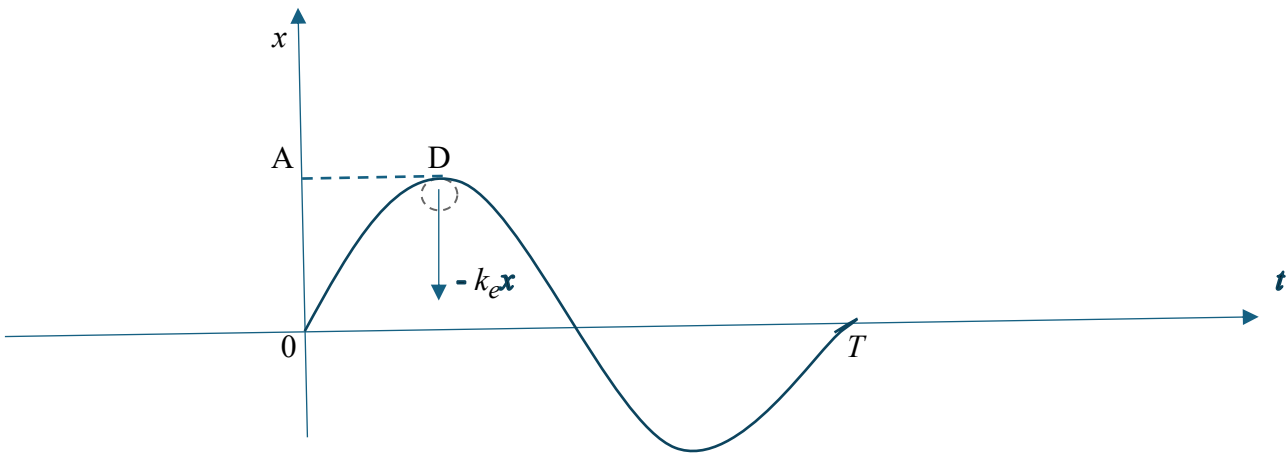

Figure 2. Social Elastic Oscillator

The aforementioned is illustrated in Fig. 2, representing the dependance of the assessment $x$ on time $t$. The designer "creative entity" is assessed through its "work", and the assessment is a circle located at the point D at a certain moment of time, when the assessment's amplitude is equal to *A* and is reached at the moment $t=T/4= \pi/2\omega=1/4\nu$. Here *T* is the period, $\nu$ is the frequency, and $\omega$ is the circular frequency.

The curve in Fig. 2 is a solution (16) of the equation (13) with $\omega = \sqrt{k_e/m}$. It is seen from the (17), (18) and Fig. 2 that $x(0)=B=0$ and $v(0)=dx/dt(0)=A\cdot\omega$. It means that the customer "others" assessment has its initial speed $v(0)$ at the moment $t=0$, e.g., right after a fashion show. When the time goes by, the force of social elasticity, opposing this fashion change, increases reaching its maximum at the moment of time $t=T/4$, and the maximum value of the assessment $x=A$ goes down after this moment.

Fig. 2 is, of course, an idealized picture of what is going on in a real community. However, it is very useful because almost any periodic function obtained from observations can be presented as a mere sum of sinusoidal functions through its expansion into a so-called Fourier series.

*Political Oscillations*

Oscillations of government power between two major independent political parties in democratic societies are well seen with a naked eye. Examples of the parties are conservative and liberal, republican and democratic, right-wing and left-wing party alliances, etc. This bipolar social structure is a natural image of the already quoted here Hegel philosophic principle of the unity and struggle of opposites. It is applicable to both alive and inanimate nature. A multipolar world concept of independent of each other countries some are trying to promote now does not look so natural and easily enforceable even if there were a strong third party in the future, such as the U.N., possessing much greater legitimacy and power to manage the entire world with no nuclear weapons at countries' disposal.

Nowadays, a multi-party semi-artificial arrangement is employed sometimes by a one-party administrative system, a disguised autocracy, to mask their political monopoly and pretend to be a fair democratic society. However, multi-party systems encounter in real democratic societies too, where different parties form alliances representing eventually the two opposites, carrying on an ongoing struggle between each other for power. The alliance compositions change making them even more adaptable to the changing world compared to their bipartisan counterparts.

To interpret the political oscillations in the framework of the presented model, it is required to specify respective transaction triangles first. In fact, they have already been identified here in the subsection *Forces* while considering examples of the Ukrainian "Maidan Revolution" of 2014 and the presidential elections in Belarus of 2020. In both cases, the elements "work" made by the government and opposition "creative entities" were respectively an autocracy and democracy.

In a case, for example, of republican-democratic oscillations in the U.S. Congress, the "others" are all U.S. voters, the "creative entity" is, respectively, republicans or democrats, while the "work" is, basically, the same – the democracy – because both parties share the same fundamental values. Differences between their "work" are minor compared to the "work" differences in the above examples of Belarus and Ukraine, earlier mentioned examples of the Red and White armies during the Russian Civil War of 1918-20, or between the North and South during the U.S. Civil War of 1861-65.

The democratic party "work" can be conditionally and very roughly called, for example, a "social democracy", addressing, first of all, average and below average Americans, while the republican "work" could be approximately defined as an "economic democracy", aimed at achieving economic freedom and prosperity through professionally applying stimulating economy measures to boost its development. Achieving this aim is also beneficial for both average and below average Americans, as well as for the entire society.

An application of the presented model will be discussed here down below in more detail, accounting for these relatively small, non-antagonistic differences. A similar sort of difference can be observed within unionized corporations between their management and union leaderships (Hussein, 2020).

The voter "others" analyze the "work" in both democratic and republican triangles and judge, in essence, through evaluating the relation (1) explicitly or implicitly by expected future benefits and costs of both parties' "work", while deciding on who they should vote for or whose "work" they should buy. This means that the situation is very similar to the example of fashion oscillations described in the previous section. While two "opposites" here are two opposing parties, it is not necessary, though is still desirable, to research on both. If the voters' assessments of the republican party satisfy the same equation (13); oscillate about like in Fig. 2.; and are correct; then their assessments of the democratic party may well be about in antiphase.

Two opposite fashion designers' "creative entity" in the described earlier example of fashion fluctuations are short-skirt and long-skirt or, respectively, narrow and wide trousers designers, but only one of them, e.g., short-skirt or narrow-trousers, was actually considered in the subsection *Fashion Fluctuations*. Because opposite parties' "work" is somewhat similar in democratic societies, the voters' choice may often be prescribed to just their changing moods towards the parties' prominent figures, expressed in their reevaluation of the leaders' ethical and/or professional dimensions. Then the respective oscillations may become described by the equation (15) instead of the (13).

*Social Sustainability*

Oscillations in business, economic, and political areas described by the second order differential equations (13) and (15) are simple and small linear motions. However, these motions are very important means to maintain economic and political stability because they relieve tensions without breaking the respective economic and political structures.

The situation is similar again to some phenomena in inanimate nature, the brightest of which is, most likely, earthquakes. Due to continuous slow motions and interactions of tectonic plates inside the Earth, there are spots in the crust, where internal stresses reach critical values, making respective elastic structures capable of breaking.

If such a structure is not very homogeneous and can be broken in many nearby places at smaller critical values, the accumulated stresses are released by microearthquakes through the so-called microseisms, and a large, destructive earthquake may never happen in this particular geographic location. The microseisms themselves are small elastic oscillations and waves described in a way similar to the presented earlier fashion and political fluctuations.

If the crust's elastic structure is more unified and a higher critical value is required to have it broken, it may well result in a single catastrophic event of releasing the accumulated high stresses, i.e., in a strong earthquake destroying the crust's structure itself and, if it happens under the ocean bottom, causing tsunamis devastating vast coastal areas (e.g., Zhdanov, 1993). These complete destructions of large material substances remind of catastrophic social phenomena, such as the so-called revolutions, causing breakages of legal administrative power along with the respective chaos and lawlessness.

Indeed, large "earthquakes" happen in human societies as well (e.g., French revolutions of 1789, 1871; Russian revolutions of 1917, 1991). However, some societies, such as England, have been gradually developing protective measures for a very long time to release people tensions in small portions, like human "microseisms". These measures include, for example, spreading the government power "horizontally" to somewhat like what is now called government branches or separation of powers, as well as "vertically" to regional levels, and, also, elaborating reliable mechanisms of power succession.

An example of developing these measures in the XVI century is queen Elizabeth the I, the last of Tudors, who was one of the greatest English monarchs. At about the same time in Russia, Ivan the Terrible, the last of Rurikovich dynasty, managed to make his power almost absolute through promoting internal terror directed to his fellow countrypeople. A result was famous long-term turmoil in Russia soon after his death until Romanov dynasty came to power.

Near present time, in the course of the last century, was even worse: a huge "earthquake", the so-called "Great October Socialist Revolution" of the 1917 along with the following it for almost 70 years socialistic dictatorship, would, most likely, scare to death even notorious medieval rulers Ivan the Terrible in Russia and Henry the VIII in England. In spite that the end of the socialism marked a turmoil of the 90th in Russia, this hopeless, oppressive administrative system is still on, unfortunately, in some societies.

Present day democratic countries utilize several well-developed mechanisms to release tensions, accumulated in their societies, in small "microseismical" portions through multiple channels, e.g., through really independent different branches of power, independent elected regional authorities, freedom of media, assembly, speech, etc. Therefore, democratic societies are much more stable than autocratic regimes, in spite that the last ones often pretend to the full unity amongst their individuals, reminding us about dangerously robust, homogeneous seismically active areas inside the Earth's crust.

Moreover, the democratic administrative system is the only fit to the market economy for the following two reasons. Firstly, if a political or economic system is somewhat shifted from its equilibrium state, then its movements – free oscillations – are governed by the same equations (13) or (15) both for the political and economic systems. I.e., free market economic and democratic political systems are both capable of releasing tensions in small portions and are more stable than other political-economic combinations.

An example of one of these other combinations is a socialistic administrative system trying to enforce a market economy. It may work well for some time, especially if it is supported from outside by a successful democratic free market society, but its either political or economic structure will eventually crumble to make a better fit to each other because it will also be a more stable combination. The socialistic political system along with the so-called planned economy adapt, obviously, better to each other. However, this society structure is less competitive than the democratic free market one. It has been persuasively shown by the breakdown of the socialistic Russia, the former so-called U.S.S.R.

Secondly, the term "free market economy" means here that exchange processes in any economic transaction triangle are free and voluntary. For example, no one "creative entity" or its criminal "representative" can force another "creative entity" to drop or raise a price for their "work" because, for example, the current price makes impossible their survival itself. Also, no social body inside a country, except for its legal political government, can demand an entity to pay them fees, in essence, additional taxes for their protection or other imposed "services" (it was called "roof" or "krysha" in the 90th Russia). Illegal autocratic governments act often exactly like organized criminal communities due to the lack of a real, not fake, opposition and, therefore, independent government branches. These actions make their free-market economy impossible.

A major, if not the major, task of all governments, including their executive, judicial, and legislative branches, is to safeguard functioning of their economies because the sense and the nature of origin of a social body's political system is to ensure sustainability of its economic development, and this is impossible without a political system in place. Therefore, a somewhat popular now expression of "politicizing the economy" does not make sense because any economy is always combined with a political system ensuring its proper functioning. There is even an area of science called Political Economy founded yet by Adam Smith.

The term "democratic political system" is equivalent here to the term "free market political system" or just "market politics" because, similar to the free market economy, it also means that political exchange processes are free and voluntary in any political transaction triangle, exactly like economic processes are free and voluntary in any business transaction triangle of a market economy. This includes the existence of more than one really competing political parties, i.e., the existence of real, not fake, opposition. This also includes basic civilian freedoms such as freedom of speech everywhere, freedom of assembly, media, and real, not fake, independence of all three government branches, etc. If a government, pretending to be a democratic one, abuses, for example, its enforcement power to keep staying in power, it won't be able to resist temptations of using coercive economic practices as well because, as it is said, "the hardest thing of all is to overpower yourself".

The following statement results from the above reasoning: existence of market politics is the necessary and sufficient condition for long-lasting existence of market economy. This seems to be true despite the well-known exclusions such as the present-day socialistic China itself and the New Economic Policy of 1921-1927 in then socialistic Russia. The latter one lasted for a very short span compared to the historic timescale and was replaced by the so-called planned economy.

If the above statement is really true, then the Chinese either political or economic system will eventually change to adopt better to each other. More importantly, the definitions themselves of capitalism and socialism accepted today must be adjusted. Indeed, it is considered now that capitalism and market economy are synonyms. Due to the above, capitalism is the market economy plus market politics, i.e., democratic administrative system, while socialism is the planned economy plus the autocratic administrative structure.

Owing to the above-described reflections, the nature of the mysteriously huge politico-economic development of the U.S. is becoming clearer. Indeed, one of the youngest countries in the world turned into the greatest economic power on Earth in less than 200 years from its emergence. The so-called American Dream has become true owing to the U.S. ingenious founding fathers, who invented and implemented a democratic political system, the market politics, blended harmonically with the best-known system for economic development and personal actualization, the market economy. The following generations have been contributing to this initial breakthrough, turning the country into an example to follow and the informal leader of the so-called free world.

CONCLUSION

The presented paper introduces subjects of the research, the space where these subjects operate, forces acting on them and affecting their motions in this space. According to the paper's name "dynamics", it proposes equations relating the forces to results of their actions – the subjects' motions – in the introduced space.

The entire construction represents a model of the subjects' variability drawn on a concept of market power originally introduced for technology commercialization (Danov, Smith, and Mitchell, 2003), which is modified and generalized to become a driving force, applicable to a much wider range of business processes, as well as to other social processes, including political ones. The proposed equations describe the effects of applying the driving forces to social bodies being participants of business and other social relationships. Therefore, the model describes relationship dynamics between participating social bodies or, in essence, dynamics of social bodies.

In addition to the concept of market power, the model is drawn on the stakeholder theory (Freeman, 1984; Mitchell, Agle, and Wood, 1997), but effects of legitimacy and urgency are not taken into account, and just effects of power are considered. Another very important "ancestor" to the presented model is the transaction cognition theory (Mitchell, 2001, 2002). Long-term continuous transaction triangles are considered in a wider range of social areas including their original area of entrepreneurship. A substantial approximation here is that just isolated triangles are examined, and potential interactions amongst triangles are briefly discussed in the framework of social damping, friction, or attenuation only.

One more foundation for this model is relationship marketing research including ethical relational constructs such as trust, commitment, communication (e.g., Anderson and Weitz, 1989; Crosby, et. al., 1990; Morgan and Hunt, 1994; Smith and Barclay, 1997), supplemented by professional or practical and other ethical or aesthetic dimensions.

These constructs are very important because unlike positions and motions of physical bodies in a regular space, visible with the naked eye, motions of social bodies take place in an invisible, introduced here Bayesian space of assessments. The position of a social body relative to an evaluator becomes visible in this space only owing to polls conducted through representative samples amongst all social bodies contained in the evaluator, being a complex social body. Results of the respective measurements are summarized to obtain the final assessment. Any independent professional or ethical variable may represent a dimension in this space, depending on a particular research problem. Therefore, this space may be multidimensional.

A suggested here range for every single evaluation relative to a contained in the evaluator social body is in the segment between [-100, +100] percentage points, while the final interval is virtually unlimited, i.e., between ($-\infty$, $+\infty$), for large complex social bodies composed of many small ones. Therefore, a new unit of measurements is introduced, and it is proposed here to name it "leo" after Leo Tolstoy, whose philosophic deliberations were a major motivation for the presented paper.

The proposed model is initial and represents just the beginnings of the new direction, deterministic social dynamics, in part because it is very hypothetical. Unlike with physical sciences, where a major statement is proven by observations that gravitational forces are indeed central, a similar statement here is a mere assumption that forces of attraction amongst social bodies are also central. While all utilized assumptions and hypotheses seem to be reasonable and relevant to our everyday experience, no field experiments have been carried out yet in the introduced here Bayesian space of assessments to quantitatively confirm or deny them.

In spite the proposed equations (8) are full nonlinear ones to describe effects of social forcing and resulting "people motions" in the Bayesian space of assessments, two specific examples of business and political processes are analyzed here in linear approximations. This linearization is proven to work well for many frequently encountered mechanical oscillations, and it should be the case here too.

It is assumed that social bodies' actions are due mostly to their feelings and moods rather than to their rational analysis, in spite that emotional and rational are interdependent (e.g., Robbins, and Judge, 2019). This seems to be consistent with human experience. Napoleon, for example, never planned on invading Russia; his major rational goal was Britain, and it was an "unexplainable chain of events", which resulted in his eventual invasion of Russia (Tolstoy, 2010).

It is also implicitly assumed that social bodies' actions follow immediately by respective changes in their relationships, which are expressed in changed assessments. While this may well be the case in real life, there also may be a lag between them and even full absence of a respective action. In addition, the same changes in judgements may well result in different degrees of intensity of the respective actions. These topics have not been fully considered here, except for some general deliberations on inertia, sustainability, or stability of social bodies.

Further experimental research could be conducted on effects of changed judgements on social bodies' actions to understand better the nature of the respective people's motions. The subject of future experimental studies could also become fluctuations in assessments of and by specific social bodies in business, political, and other functional areas. Attenuation may also be considered to bring the model results closer to respective observations. The assessments could be measured in "leos" through using representative samples to conduct the respective polls.

Further theoretical research may be devoted to verifying the introduced here hypothesis that market politics is the necessary and sufficient condition for market economy as well as to identify the nature and duration of possible exclusions from its applicability.

Research on measuring masses of social bodies would enable one to estimate coefficients of attraction and doubts or social elasticity. This, in turn, would make it possible to evaluate social forces of attraction and doubts/elasticity as well as magnitudes of other forces applied to social bodies. Finally, nonlinear effects could be specifically studied to understand better some destructive social phenomena such as insurrections, wars, revolutions, etc.

Altogether, successful future research in these areas may well result in the emergence of new quantitative modelling of business, political, and other social processes. Drawn on the presented here idealized social dynamics, this development of numerical modelling methods for social phenomena, including destructive ones, can take social forecasting to another, much higher level, compared to what it is now.

## ACKNOWLEDGEMENTS

Thank you to Professors Ronald K. Mitchell and Sergey E. Mikhailov for their comments which helped to improve the paper.